# Sum Rate Optimization for Coordinated Multi-Antenna Base Station Systems

Tadilo Endeshaw and Luc Vandendorpe
ICTEAM Institute
Université catholique de Louvain
Place du Levant, 2, B-1348, Louvain La Neuve, Belgium
Email: {tadilo.bogale, luc.vandendorpe}@uclouvain.be

*Abstract*— This paper considers the joint precoder design problem for multiple-input single-output (MISO) systems with coordinated base stations (BSs). We consider maximization of the total sum rate with per BS antenna power constraint problem. For this problem, we propose a novel linear iterative algorithm. The problem is solved as follows. First, by introducing additional optimization variables and applying matrix fractional minimization, we reformulate the original problem as a new problem. Second, for the given precoder vectors of all users, we optimize the introduced variables of the latter problem using Geometric Programming (GP) method. Third, keeping the introduced variables constant, the precoder vectors of all users are optimized by using phase rotation technique. The second and third steps are repeated until convergence. We have shown that the proposed algorithm is guaranteed to converge. Moreover, for the total sum power constraint case, simulation results show that the proposed iterative algorithm achieves almost the same average performance as that of the algorithm which utilizes mean-square-error (MSE) uplink-downlink duality approach. We also show that our iterative algorithm can be used to solve sum rate maximization and weighted sum MSE minimization problems for an arbitrary power constraint.

## I. Introduction

Multi-antenna systems have been proven to enhance the spectral efficiency of wireless systems. This performance improvement is achieved by employing signal processing at the transmitters (precoder) and receivers (decoders).

The achievable sum rate of the broadcast channel (BC) obtained by dirty paper precoding technique has been characterized for multiple-input single-output (MISO) systems [1]. The latter work has been extended in [2] for multiple-input multiple-output (MIMO) systems. The authors of [3] and [4] have shown that "dirty paper coding" (DPC) achieves the capacity region of BC channels. However, due to the non-linear characteristics of DPC, practical realization of DPC has appeared to be difficult.

Given the drawbacks of DPC, linear processing is motivated as it exhibits good performance/complexity trade-off. However, finding linear processing schemes that achieve the capacity of BC channels is still an open issue. In [5], linear processing method such as channel block-diagonalization is suggested. This method suffers from noise enhancement and has a restriction on the number of transmit and receive antennas. In [6], weighted sum rate maximization problem is formulated as the problem of minimizing the geometric product of minimum mean-square-errors (GPMMSE). To solve the optimization problem an iterative approach which uses mean-square-error (MSE) uplink-downlink duality is suggested. Minimizing the product of all users minimum mean-square-error (MMSE) matrix determinants is proposed as an equivalent formulation for the un-weighted sum rate maximization [7]. This problem is non-convex and it is solved by employing sequential quadratic programming (SQP). In [8], the robust sum rate maximization problem has been examined. The latter paper examines the problem using the worst-case robust design approach, and utilizes MSE uplink-downlink duality approach to solve the sum rate maximization problem.

All of the above papers examine their problems for conventional cellular networks with a total base station (BS) power constraint. In these networks, BSs from different cells communicate with their respective remote terminals independently. Hence, inter-cell interference is considered as a background noise. Recently, it has been shown that BS-cooperative communication is a promising technique to mitigate inter-cell interference [9], [10], [11], [12]. In [9], four MSE-based linear transceiver optimization problems have been considered for the coordinated BS MIMO systems. These problems are examined by assuming that the total power of each BS or the individual power of each BS antenna (group of antennas) is constrained. The optimization problems in [9] are solved as follows: first, by keeping the receivers constant, optimization of the precoder matrices are formulated as a second order cone program (SOCP) problem (SOCP problems are convex and can be solved by using existing convex optimization tools). Second, for the given BS precoders, the receiver of each user is optimized by MMSE technique. These steps are repeated iteratively to jointly optimize the transmitters and receivers.

In this paper, we examine the joint optimization of the precoders for maximizing the total sum rate with per BS antenna power constraint[1]. To solve this problem, first, by intro-

The authors would like to thank the Region Wallonne for the financial support of the project MIMOCOM in the framework of which this work has been achieved.

---

[1]According to [13], in a practical multi-antenna BS systems, each BS antenna has its own power amplifier and the maximum power for each BS antenna is limited by some value. This motivates us to consider the power constraint of each BS antenna. However, as will be clear later, our proposed algorithm can be modified straightforwardly to handle any arbitrary power constraint.

ducing additional optimization variables and applying matrix fractional minimization, we reformulate the original problem as a new problem. Second, for the given precoder vectors of all users, we optimize the introduced variables of the latter problem using Geometric Programming (GP) method. Third, keeping the introduced variables constant, the precoder vectors of all users are optimized by using phase rotation technique. The second and third steps are repeated until convergence. We have shown that the proposed algorithm is guaranteed to converge. Moreover, for the total sum power constraint case, simulation results show that the proposed iterative algorithm achieves almost the same average performance as that of the algorithm which utilizes MSE uplink-downlink duality approach (see [6] and [8] for the details of this approach). We have also shown that our iterative algorithm can be used to solve weighted sum MSE minimization problem. In this paper, it is assumed that each BS is equipped with multiple antennas and the users are assumed to have single antenna, and perfect channel state information (CSI) is available both at the BSs and mobile stations (MSs).

## II. SYSTEM MODEL

We consider a coordinated BS system where $L$ BSs are serving $K$ decentralized single antenna MSs. The $l$th BS is equipped with $N_l$ transmit antennas. By denoting the symbol intended for the $k$th user as $d_k$, the entire symbol can be written in a data vector $\mathbf{d} \in \mathcal{C}^{K \times 1}$ as $\mathbf{d} = [d_1, \cdots, d_K]^T$. The $l$th BS precodes $\mathbf{d}$ into an $N_l$ length vector by using its overall precoder matrix $\mathbf{B}_l = [\mathbf{b}_{l1}, \cdots, \mathbf{b}_{lK}]$, where $\mathbf{b}_{lk} \in C^{N_l \times 1}$ is the precoder vector of the $l$th BS for the $k$th MS. The $k$th MS employs a receiver $w_k$ to estimate its symbol $d_k$. The estimated symbol at the $k$th MS is given by

$$\hat{d}_k = w_k^H (\mathbf{h}_k^H \sum_{i=1}^{K} \mathbf{b}_i d_i + n_k) \quad (1)$$

where $\mathbf{h}_k^H = [\mathbf{h}_{1k}^H, \cdots, \mathbf{h}_{Lk}^H] \in \mathcal{C}^{1 \times N}$ with $\mathbf{h}_{lk}^H \in \mathcal{C}^{1 \times N_l}$ as the channel vector between the $l$th BS and the $k$th MS, $\mathbf{b}_k = [\mathbf{b}_{1k}^T, \cdots, \mathbf{b}_{Lk}^T]^T \in \mathcal{C}^{N \times 1}$, $N = \sum_{l=1}^{L} N_l$ and $n_k$ is the additive noise at the $k$th MS. It is clearly seen that the expression (1) has exactly the same form as the estimate of $d_k$ for the downlink MISO system where a BS equipped with $N$ transmitt antennas is serving $K$ decentralized users. Hence, we can interpret the coordinated BS system as a one giant downlink system [9], [12]. For convenience, we follow the same channel vector notations as in [14]. It is assumed that $n_k$ is a zero-mean circularly symmetric complex Gaussian (ZMCSCG) random variable with variance $\sigma_k^2$, i.e., $n_k \sim \mathcal{NC}(0, \sigma_k^2)$. We also assume that the symbol $d_k$ is a ZMCSCG random variable with unit variance and is independent of $\{d_i\}_{i=1, i\neq k}^{K}$ and noise $n_k$, i.e., $\mathrm{E}\{d_k d_k^H\} = 1$, $\mathrm{E}\{d_k d_i^H\} = 0$, $\forall i \neq k$ and $\mathrm{E}\{d_k n_k^H\} = 0$. For this system model, when perfect CSI is available at the BS and MSs, the $k$th user achievable rate is given by

$$R_k = \log_2(1 + \psi_k) \quad (2)$$

where $\psi_k = \frac{\mathbf{h}_k^H \mathbf{b}_k \mathbf{b}_k^H \mathbf{h}_k}{\mathbf{h}_k^H \sum_{i=1, i\neq k}^{K} \mathbf{b}_i \mathbf{b}_i^H \mathbf{h}_k + \sigma_k^2}$.

## III. PROBLEM FORMULATION AND PROPOSED SOLUTION

Mathematically, the sum rate maximization problem can be formulated as

$$\max_{\{\mathbf{b}_k\}_{k=1}^{K}} \sum_{k=1}^{K} \log_2(1 + \psi_k), \text{ s.t } [\sum_{k=1}^{K} \mathbf{b}_k \mathbf{b}_k^H]_{n,n} \le p_n, \ \forall n \quad (3)$$

where $p_n$ is the maximum power allocated to the $n$th antenna of all BSs. After straightforward mathematical manipulations, the latter problem can be equivalently expressed as

$$\min_{\{\mathbf{b}_k\}_{k=1}^{K}} \prod_{k=1}^{K}(1 + \psi_k)^{-1}, \text{ s.t } [\sum_{k=1}^{K} \mathbf{b}_k \mathbf{b}_k^H]_{n,n} \le p_n, \ \forall n. \quad (4)$$

The above optimization problem is not convex. Thus, convex optimization techniques can not be applied to (4). To this end, we propose an iterative approach to solve this problem. In this regard, we consider the following Lemma.

*Lemma 1*: The suboptimal $\{\mathbf{b}_k\}_{k=1}^{K}$ of (4) can be obtained by solving the following problem.

$$\min_{\{\mathbf{b}_k, \nu_k\}_{k=1}^{K}} \left( \frac{1}{K} \sum_{k=1}^{K} \nu_k (1 + \psi_k)^{-1} \right)^{K}$$

$$\text{s.t } [\sum_{k=1}^{K} \mathbf{b}_k \mathbf{b}_k^H]_{n,n} \le p_n, \ \prod_{k=1}^{K} \nu_k = 1, \nu_k \ge 0, \ \forall n, k. \quad (5)$$

In particular, for $K = 2$, this problem can be equivalently formulated by

$$\min_{\nu, \{\mathbf{b}_k\}_{k=1}^{2}} \frac{1}{2} \left( \nu[(1 + \psi_1)^{-1}]^2 + \frac{1}{\nu}[(1 + \psi_2)^{-1}]^2 \right)$$

$$\text{s.t } [\sum_{k=1}^{K} \mathbf{b}_k \mathbf{b}_k^H]_{n,n} \le p_n, \ \forall n, \ \nu \ge 0. \quad (6)$$

*Proof 1*: For fixed $\{\mathbf{b}_k\}_{k=1}^{K}$, optimizing $\{\nu_k\}_{k=1}^{K}$ of (5) can be expressed as

$$\min_{\{\nu_k\}_{k=1}^{K}} \left( \frac{1}{K} \sum_{k=1}^{K} \beta_k \nu_k \right)^{K}$$

$$\text{s.t } \prod_{k=1}^{K} \nu_k = 1, \nu_k \ge 0, \ \forall k. \quad (7)$$

where $\beta_k = (1 + \psi_k)^{-1}$. The above problem is GP problem for which global optimality is guaranteed. Clearly, the optimal solution of (7) satisfy $\{\nu_k > 0\}_{k=1}^{K}$, and the objective and constraint functions of this problem are continuously differentiable. Moreover, by replacing $\nu_1 = (\prod_{k=2}^{K} \nu_k)^{-1}$, the equality constraint of the latter problem can be removed. These two facts show that the optimal $\{\nu_k\}_{k=1}^{K}$ of the above problem are regular [15], [16]. Thus, the global optimal solution of (7) can be obtained by choosing $\{\nu_k\}_{k=1}^{K}$ that satisfy the Karush-Kuhn-Tucker (KKT) optimality conditions which are given by [17]

$$\beta_k \left( \frac{1}{K} \sum_{i=1}^{K} \beta_i \nu_i \right)^{K-1} - \gamma \prod_{i=1, i\neq k}^{K} \nu_i - \lambda_k = 0 \quad (8)$$

$$\lambda_k \nu_k = 0 \quad (9)$$

$$\lambda_k \ge 0 \quad (10)$$

where $\gamma$ and $\{\lambda_k\}_{k=1}^K$ are the Lagrangian multipliers corresponding to the constraints $\prod_{k=1}^K \nu_k = 1$ and $\{\nu_k \geq 0\}_{k=1}^K$, respectively. Now, multiplying (8) by $\nu_k$ and employing $\prod_{k=1}^K \nu_k = 1$ and (9), we get

$$\beta_k \nu_k \left(\frac{1}{K}\sum_{i=1}^K \beta_i \nu_i\right)^{K-1} - \gamma \nu_k \prod_{i=1, i\neq k}^K \nu_i - \lambda_k \nu_k = 0$$
$$\Rightarrow \beta_k \nu_k \left(\frac{1}{K}\sum_{i=1}^K \beta_i \nu_i\right)^{K-1} = \gamma \prod_{i=1}^K \nu_i = \gamma. \quad (11)$$

By summing the K equalities of (11), $\gamma$ can be determined by

$$\gamma = \left(\frac{1}{K}\sum_{i=1}^K \beta_i \nu_i\right)^K. \quad (12)$$

Substituting $\gamma$ of (12) into (11), and noting that $\frac{1}{K}\sum_{k=1}^K \beta_k \nu_k > 0$ we obtain

$$\beta_k \nu_k = \frac{1}{K}\sum_{i=1}^K \beta_i \nu_i. \quad (13)$$

Multiplying the $K$ equalities of (13) yields

$$\prod_{k=1}^K \beta_k = \left(\frac{1}{K}\sum_{i=1}^K \beta_i \nu_i\right)^K$$
$$\Rightarrow \prod_{k=1}^K (1+\psi_k)^{-1} = \left(\frac{1}{K}\sum_{k=1}^K \nu_k (1+\psi_k)^{-1}\right)^K. \quad (14)$$

The above expression shows that the suboptimal solution of (4) can be equivalently obtained by solving (5). For $K=2$, by employing (13) and (14), we can express (5) as in (6)[18].

Since the solution of (5) satisfies $\frac{1}{K}\sum_{k=1}^K \nu_k(1+\psi_k)^{-1} > 0$, the objective function of the latter problem can be replaced by $\sum_{k=1}^K \nu_k(1+\psi_k)^{-1}$. As a result, (5) can be equivalently expressed as

$$\min_{\{\mathbf{b}_k,\nu_k\}_{k=1}^K} \sum_{k=1}^K \nu_k(1+\psi_k)^{-1}$$
$$\text{s.t } [\sum_{k=1}^K \mathbf{b}_k \mathbf{b}_k^H]_{n,n} \leq p_n, \quad \prod_{k=1}^K \nu_k = 1, \nu_k \geq 0, \quad \forall n,k. \quad (15)$$

This problem is still not convenient to solve. Towards this end, we consider the following Lemma.

*Lemma 2*: The suboptimal $\{\mathbf{b}_k,\nu_k\}_{k=1}^K$ of (15) can be obtained by solving the following problem.

$$\min_{\{\mathbf{b}_k,\nu_k,t_k,z_k\}_{k=1}^K} \sum_{k=1}^K \nu_k \left(\frac{t_k^H t_k (\mathbf{h}_k^H \sum_{i=1,i\neq k}^K \mathbf{b}_i \mathbf{b}_i^H \mathbf{h}_k + \sigma_k^2)}{\mathbf{h}_k^H \mathbf{b}_k \mathbf{b}_k^H \mathbf{h}_k} + z_k^H z_k\right)$$
$$\text{s.t } [\sum_{k=1}^K \mathbf{b}_k \mathbf{b}_k^H]_{n,n} \leq p_n, \quad \forall n$$
$$t_k + z_k = 1, \quad \prod_{k=1}^K \nu_k = 1, \quad \forall k. \quad (16)$$

*Proof*: For given $\{\mathbf{b}_k, \nu_k\}_{k=1}^K$, the optimal $\{t_k, z_k\}_{k=1}^K$ of the above problem can be obtained by solving the following problem

$$\min_{\{t_k,z_k\}_{k=1}^K} \sum_{k=1}^K \nu_k \left(\frac{t_k^H t_k}{\psi_k} + z_k^H z_k\right), \quad \text{s.t } t_k + z_k = 1, \forall k. \quad (17)$$

Since the above problem is convex, it can be solved by using the Lagrangian multiplier method. The Lagrangian function associated with (17) is given by

$$L = \sum_{i=1}^K \left[\nu_i\left(\frac{t_i^H t_i}{\psi_i} + z_i^H z_i\right) - \mu_i(t_i + z_i - 1)\right] \quad (18)$$

where $\{\mu_i\}_{i=1}^K$ are the Lagrangian multipliers corresponding to the constraints $\{t_i + z_i = 1\}_{i=1}^K$. Taking the first order derivative of $L$ with respect to $t_k$ and $z_k$ results

$$t_k = \frac{\psi_k \mu_k}{\nu_k}, \quad z_k = \frac{\mu_k}{\nu_k}, \quad \forall k. \quad (19)$$

By substituting the above expression in the constraint of (17), we get

$$t_k + z_k = 1 \Rightarrow \frac{\psi_k \mu_k}{\nu_k} + \frac{\mu_k}{\nu_k} = 1 \Rightarrow \mu_k = \frac{\nu_k}{1+\psi_k}.$$

It follows that

$$t_k = \frac{\psi_k}{1+\psi_k}, \quad z_k = \frac{1}{1+\psi_k}, \quad \forall k. \quad (20)$$

Plugging (20) into the objective function of (17) yields the objective function of (15). It follows that the optimal solution of (15) can be obtained by solving (16). Lemma 2 is proved by modifying the technique of matrix fractional minimization in [17] and [19].

**Note**: The optimal solution of (15) might satisfy $\mathbf{h}_k^H \mathbf{b}_k \mathbf{b}_k^H \mathbf{h}_k = 0, \exists k$, (which means that these users are switched off). In (16), the terms $\{\mathbf{h}_k^H \mathbf{b}_k \mathbf{b}_k^H \mathbf{h}_k\}_{k=1}^K$ are not allowed to be zero, However, the latter terms can be arbitrarily very close to zero. In practice, if $\mathbf{h}_k^H \mathbf{b}_k \mathbf{b}_k^H \mathbf{h}_k$ is below a certain threshold, the corresponding user can be considered as switched off. This shows that the solution obtained by (16) can identify switched off users. Nonetheless, reformulating (15) by considering $\{\mathbf{h}_k^H \mathbf{b}_k \mathbf{b}_k^H \mathbf{h}_k \geq 0\}_{k=1}^K$ is still an open problem.

Again (16) can not be solved in its current form. To simplify this problem, we apply Lemma 1 for $K=2$ case twice by introducing $\{\tau_k\}_{k=1}^K$ and $\{\eta_k\}_{k=1}^K$, respectively. Upon doing so, (16) can be reexpressed as

$$\min_{\{\mathbf{b}_k,t_k,\nu_k,\tau_k,\eta_k\}_{k=1}^K} \sum_{k=1}^K \nu_k \left[(t_k^H-1)(t_k-1) + \frac{\tau_k}{2(\mathbf{h}_k^H \mathbf{b}_k \mathbf{b}_k^H \mathbf{h}_k)^2} + \frac{1}{2\tau_k}\left(\frac{(t_k^H t_k)^4}{2\eta_k} + \frac{\eta_k}{2}(\mathbf{h}_k^H \sum_{i=1,i\neq k}^K \mathbf{b}_i \mathbf{b}_i^H \mathbf{h}_k + \sigma_k^2)^4\right)\right]$$
$$\text{s.t } [\sum_{k=1}^K \mathbf{b}_k \mathbf{b}_k^H]_{n,n} \leq p_n, \quad \prod_{k=1}^K \nu_k = 1, \quad \forall n,k. \quad (21)$$

For fixed $\{\mathbf{b}_k, t_k\}_{k=1}^{K}$, the above optimization problem is a GP and thus it can be solved using standard convex optimization tool [20]. Next, keeping $\{\nu_k, \tau_k, \eta_k\}_{k=1}^{K}$ constant, the optimal $\{\mathbf{b}_k, t_k\}_{k=1}^{K}$ of (21) can be obtained as follows. For any $\{\theta_k\}_{k=1}^{K}$, we have $\mathbf{h}_k^H \mathbf{b}_k \mathbf{b}_k^H \mathbf{h}_k = |\mathbf{h}_k^H \mathbf{b}_k|^2 = |\mathbf{h}_k^H \mathbf{b}_k e^{j\theta_k}|^2$. Now, without loss of generality, by choosing $\{\theta_k\}_{k=1}^{K}$ such that $\mathbf{h}_k^H \mathbf{b}_k > 0$,[2] we can express $\mathbf{h}_k^H \mathbf{b}_k \mathbf{b}_k^H \mathbf{h}_k = (\mathbf{h}_k^H \mathbf{b}_k)^2$. By doing so and after some mathematical manipulations, the optimal $\{\mathbf{b}_k, t_k\}_{k=1}^{K}$ of (21) can be obtained by solving the following problem

$$\min_{\{\mathbf{b}_k, t_k\}_{k=1}^{K}} \sum_{k=1}^{K} \nu_k \left[ \frac{\tau_k}{2 c_k^4} + \frac{1}{2\tau_k} \left( \frac{f_k^4}{2\eta_k} + \frac{\eta_k}{2} r_k^4 \right) + x_k \right]$$

$$\text{s.t } [\sum_{k=1}^{K} \mathbf{b}_k \mathbf{b}_k^H]_{n,n} \leq p_n, \quad \forall n$$

$$(t_k^H - 1)(t_k - 1) \leq x_k, \quad t_k^H t_k \leq f_k$$

$$\mathbf{h}_k^H \sum_{i=1, i\neq k}^{K} \mathbf{b}_i \mathbf{b}_i^H \mathbf{h}_k + \sigma_k^2 \leq r_k, \quad \mathbf{h}_k^H \mathbf{b}_k = c_k, \quad \forall k. \quad (22)$$

It can be shown that for fixed $\{\nu_k, \tau_k, \eta_k\}_{k=1}^{K}$, the above problem is convex for which global optimal solution can be obtained by using interior point methods [17]. In summary, problem (21) can be solved iteratively as follows.

**Algorithm I:** Iterative algorithm to solve (21).
Initialization: Set $\{\mathbf{b}_k = \mathbf{h}_k\}_{k=1}^{K}$ and normalize $\{\mathbf{b}_k\}_{k=1}^{K}$ such that each antenna power constraint is satisfied with equality. Then, compute $\{t_k = \psi_k/(1+\psi_k)\}_{k=1}^{K}$.
**repeat**
1) With the current $\{\mathbf{b}_k, t_k\}_{k=1}^{K}$, compute the optimal $\{\nu_k, \tau_k, \eta_k\}_{k=1}^{K}$ using (21).
2) With the optimal $\{\nu_k, \tau_k, \eta_k\}_{k=1}^{K}$ obtained from the above step, calculate the optimal solution of $\{\mathbf{b}_k, t_k\}_{k=1}^{K}$ by solving (22). As we can see the solution of (21) and (22) do not let $\{\tau_k, \eta_k, \mathbf{h}_k^H \mathbf{b}_k\}_{k=1}^{K}$ to be zero. However, when one of the the terms $\{\tau_k, \eta_k, \mathbf{h}_k^H \mathbf{b}_k\}_{k=1}^{K}$ are very small (i.e., $< 10^{-6}$), we have noticed numerical instability. To handle such a problem, we replace $\{\tau_k, \eta_k, \mathbf{h}_k^H \mathbf{b}_k < 10^{-6}\}_{k=1}^{K}$ by $\{\tau_k, \eta_k, \mathbf{h}_k^H \mathbf{b}_k = 10^{-6}\}_{k=1}^{K}$.
3) Compute the objective function of (21).
**Until** Convergence.

**Convergence:-** At each step, the objective function of (21) is non-increasing and it is lower bounded by 0. These two facts show that the proposed iterative algorithm is always guaranteed to converge. It can be shown that if (21) is convergent, the original problem (3) is also convergent. The details are omitted due to space constraint. However, since the latter problem is non-convex, we are not able to show the global optimality of **Algorithm I** analytically.
**Initialization:-** In general, different initializations affect the convergence speed of **Algorithm I**. In most of our simulations, we obtain faster convergence speed when the initialization

[2] We would like to mention here that the idea of phase rotation has been used to solve different kinds of problems.(see [21] for other problems)

is performed as in **Algorithm I**. However, getting the best initialization that results the fastest convergence speed of **Algorithm I** is an open research topic.

## IV. EXTENSION TO WEIGHTED SUM MSE MINIMIZATION PROBLEM

The weighted sum MSE minimization problem can be expressed as

$$\min_{\{\mathbf{b}_k\}_{k=1}^{K}} \sum_{k=1}^{K} \upsilon_k (1+\psi_k)^{-1}, \text{ s.t } [\sum_{k=1}^{K} \mathbf{b}_k \mathbf{b}_k^H]_{n,n} \leq p_n, \ \forall n \quad (23)$$

where $\upsilon_k$ is the MSE weighting factor and $p_n$ is the power allocated to each BSs antenna. It can be shown that this problem can be solved using the same approach as the sum rate optimization problem (4). It can be clearly seen that **Algorithm I** can be modified straightforwardly to handle any arbitrary power constrained sum rate maximization (weighted sum MSE minimization) problems. We would like to mention here that for the total BS power constraint case, the authors of [22] establish the relation between the weighted sum rate maximization and weighted sum MSE minimization problems by exploiting the Lagrangian functions of these two problems. Thus, our work generalizes the relation between these two problems for any arbitrary power constraint without resorting the Lagrangian functions.

## V. SIMULATION RESULTS

In this section, we provide the simulation results for problem (4). We consider a system with $L = 2$ BSs where each BS has 2 antennas and $K = 4$ MSs. The channel between each BS and MSs consists of ZMCSCG entries with unit variance. We assume that $\{\sigma_k^2 = \sigma^2\}_{k=1}^{K}$.

The optimal transmit powers of each antenna are plotted in Fig. 1 for the case where the power constraint of each antenna is set to 2, i.e., $\{p_n = 2\}_{n=1}^{4}$ and $\sigma^2 = 0.1$. The latter figure shows that for problem (4), all antennas do not necessary utilize their full powers to maximize the total achievable rate of the system. This observation fits with that of [9] where sum MSE minimization problem has been examined for coordinated BSs with per BS power constraint.

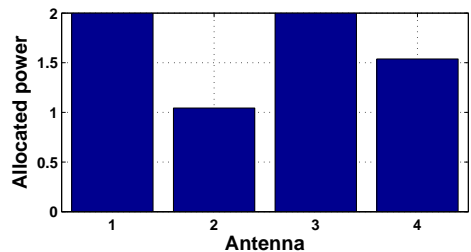

Fig. 1. The allocated power of each antenna with **Algorithm I**.

It has been shown in [6] and [8] that the sum rate maximization constrained with a total BS power problem has been solved using MSE uplink-downlink duality technique. Moreover, the simulation results of [6] has demonstrated that the duality approach of solving the sum rate optimization problem outperforms other linear schemes. This motivates us to compare the performance of [6] and our algorithm for the total sum power constraint case[3]. The comparison is performed by averaging over 1000 randomly chosen channel realizations. The signal-to-noise ratio (SNR) is defined as $P_{\text{sum}}/\sigma^2$, where $P_{\text{sum}}$ is the total sum power of all antennas and $\sigma^2$ is the noise variance. The SNR is controlled by varying $\sigma^2$ while setting $P_{\text{sum}} = 10$. Fig. 2 shows that for the total sum power constraint case, the proposed algorithm and the algorithm of [6] have almost the same average total sum rate. Note that although the algorithm in [6] and our algorithm yield very close average performance, this scenario is not always true for each channel realization and SNR value.

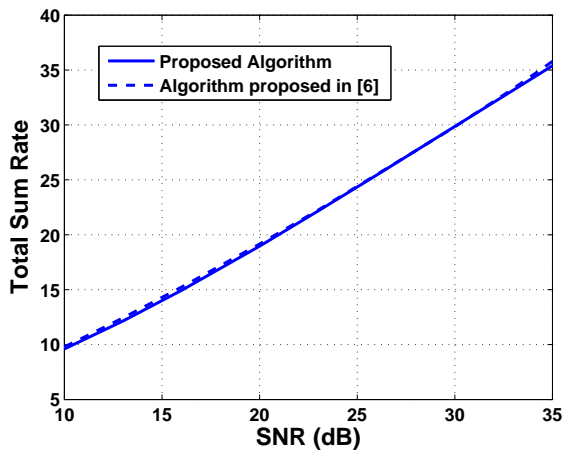

Fig. 2. Comparison of the total achievable sum rates obtained in [6] and **Algorithm I**.

## VI. CONCLUSIONS

This paper considers the precoder design problem for MISO systems with coordinated BSs. We examine maximization of the total sum rate with per BS antenna power constraint problem. This problem is efficiently solved by using our linear iterative algorithm. Our new method of precoder design employs the reformulation of the original problem into another problem, modified matrix fractional minimization, phase rotation and an iterative approach. Unlike the MSE uplink-downlink duality solution approach (this approach can solve the sum rate maximization problem only for a total BS power constraint case), our proposed approach is able to solve sum rate maximization problem with per antenna/groups of antenna BS power constraint. Moreover, for the total BSs power constraint case, simulation results have shown that the proposed iterative algorithm achieves almost the same average sum rate as that of the algorithm which utilizes the MSE uplink-downlink duality approach. We also show that our iterative algorithm can also be used to solve the weighted sum MSE minimization problem.

---

[3]We would like to mention here that the uplink-downlink duality approach of [6] can not be applied to solve (4) with per antenna BS power constraint.